\def\date{March 14, 2007}
\documentclass[12pt]{article}

\usepackage{amssymb}
\usepackage{amsmath}
\input{liemacs.sty}

\begin{document}

\newcommand{\bS}{{\mathbb S}}
\newcommand{\Fl}{{\rm Fl}}
\newcommand{\Gr}{{\rm Gr}}
\newcommand{\Map}{\mathrm{Map}}
\newcommand{\D}{\frac{d}{d\lambda}\Big|_{\lambda=0}}

\title{An implicit function theorem for Banach spaces and some applications}

\author
{Jinpeng An\footnote{ETH Z\"urich, R\"amistrasse 101, 8092 Z\"urich,
Switzerland, jinpeng.an@math.ethz.ch} and Karl-Hermann
Neeb\footnote{Technische Universit\"at Darmstadt,
Schlossgartenstrasse 7, D-64289 Darmstadt, Deutschland,
neeb@mathematik.tu-darmstadt.de}}


\maketitle

\begin{abstract}
We prove a generalized implicit function theorem for Banach spaces,
without the usual assumption that the subspaces involved being
complemented. Then we apply it to the problem of parametrization of
fibers of differentiable maps, the Lie subgroup problem for Banach--Lie
groups, as well as Weil's local rigidity for homomorphisms from
finitely generated groups to Banach--Lie groups.
\end{abstract}

\section{Introduction}

Motivated by the problem of local rigidity of homomorphisms from
finitely generated groups to (finite-dimensional) Lie groups, Weil
\cite{We64} proved the following theorem.

\begin{theorem}\label{thm:Weil0}{\rm(Weil)}
Let $L, M$ and $N$ be finite-dimensional $C^1$-manifolds,
$\varphi:L\rightarrow M$ and $\psi:M\rightarrow N$
be $C^1$-maps. Let $x\in U$, $y=\varphi(x)$. Suppose \\
{\rm(1)} $\psi\circ\varphi\equiv z\in N$;\\
{\rm(2)} $\im(d\varphi(x))=\ker(d\psi(y))$.\\
Then there exists a neighborhood $U$ of $y$ in $M$ such that
$$\psi^{-1}(z)\cap U=\varphi(L)\cap U.$$
\end{theorem}

Weil's proof relies heavily on the Implicit Function Theorem. If $L,
M$ and $N$ in Theorem~\ref{thm:Weil0} are Banach manifolds, and
$\ker(d\varphi(x)), \ker(d\psi(y))$ and $\im(d\psi(y))$ are closed
complemented subspaces in the corresponding tangent spaces of the
manifolds respectively, then Weil's proof still works,  based on the
Implicit Function Theorem for Banach spaces. Using the Nash-Moser
Inverse Function Theorem, Hamilton~\cite{Ha77} proved a similar
result in the setting of tame Fr\'{e}chet spaces under more
splitting assumptions (see also \cite{Ha82}). Such kind of results
are useful in problems of deformation rigidity (or local rigidity)
of certain geometric structures (see e.g.\
\cite{Ha82},\cite{Be00}, \cite{Fi05}, \cite{AW05}).
But such complementation or splitting
conditions are not always satisfied and hard to verify in many
concrete situations.

The first goal of this paper is to prove the Banach version of
Theorem~\ref{thm:Weil0} without the assumption that the subspaces
involved are complemented. The only additional assumption we have to
make, comparing with the finite-dimensional case, is that
$\im(d\psi(y))$ is closed. More precisely, we will prove the
following theorem. For convenience, we state its local form.

\begin{theorem}\label{thm:implicit0}
Let $X, Y$ and $Z$ be Banach spaces, $U\subeq X$, $V\subeq Y$ be
open $0$-neighborhoods. Let $\varphi:U\rightarrow V$ and
$\psi:V\rightarrow Z$ be $C^1$-maps.
Suppose\\
{\rm(1)} $\varphi(0)=0$;\\
{\rm(2)} $\psi\circ\varphi\equiv0$;\\
{\rm(3)} $\im(d\varphi(0))=\ker(d\psi(0))$;\\
{\rm(4)} $\im(d\psi(0))$ is closed.\\
Then there exists a $0$-neighborhood $W\subeq V$ such that
$$\psi^{-1}(0)\cap W=\varphi(U)\cap W.$$
\end{theorem}

Informally, we may think Theorem~\ref{thm:implicit0} as saying that
$\varphi$ determines some ``implicit function" from
$\psi^{-1}(0)\cap W$ to $U$ (although we do not even know a prori if it is
continuous). So we call Theorem~\ref{thm:implicit0} a
``generalized implicit function theorem".

We will give three applications of the generalized implicit function
theorem.

(1) \textbf{Parametrization of fibers of differentiable maps.}
Suppose $M, N$ are finite-dimensional $C^k$-manifolds and
$\psi:M\rightarrow N$ is a $C^k$-map. If $z\in N$ is a regular value
of $\psi$ in the sense that for each $y\in\psi^{-1}(z)$, the
differential $d\psi(y):T_yM\rightarrow T_zN$ is surjective, then the
Implicit Function Theorem implies that $\psi^{-1}(z)$ is a
$C^k$-submanifold of $M$. This fact can be generalized directly to
Banach manifolds under the additional assumption that
$\ker(d\psi(y))$ is complemented in $T_yM$. Without the
complementation assumption, we can not conclude that $\psi^{-1}(z)$
is a submanifold in general. However, using the generalized implicit
function theorem, we can prove that under the assumption of the
existence of the map $\phi$ as in Theorem~\ref{thm:implicit0},
$\psi^{-1}(z)$ can be parametrized as a $C^k$-manifold in such a way
that the manifold structure is compatible with the subspace
topology, and the inclusion map is $C^k$.

(2) \textbf{Lie subgroups of Banach--Lie groups.} If $G$ is a Banach
Lie group, we call a closed subgroup $H$ of $G$ a Lie subgroup if
there is a Banach--Lie group structure on $H$ compatible with the
subspace topology (there is at most one such structure).
It is well-known that closed subgroups of
finite-dimensional Lie groups are Lie subgroups but this is false for
Banach--Lie groups (cf.\ \cite{Hof75}; \cite{Ne06}, Rem.~IV.3.17). Using the
generalized implicit function theorem and a theorem of Hofmann
\cite{Hof75}, we will prove that certain isotropy groups are Lie
subgroups. In particular, we will prove that if the coset space
$G/H$ carries a Banach manifold structure for which $G$ acts
smoothly and the quotient map $q \: G \to G/H, g \mapsto gH$ has
surjective differential in some point $g\in G$, then $H$ is a
Banach--Lie subgroup of $G$.

(3) \textbf{Local rigidity of homomorphisms.} Let $G$ be a
Banach--Lie group, $\Gamma$ a finitely generated group. A
homomorphism $r:\Gamma\to G$ is locally rigid if any homomorphism
from $\Gamma$ to $G$ sufficiently close to $r$ is conjugate to $r$.
Using Theorem~\ref{thm:Weil0}, Weil \cite{We64} proved that if
$H^1(\Gamma,\L(G))=\{0\}$, then $r$ is locally rigid. We will
generalize Weil's result to the Banach setting under the assumption
that a certain linear operator has closed range.

Theorem~\ref{thm:implicit0} will be proved in Section 2. The
aforementioned three applications will be derived in Sections
3--5, respectively.

This work was initiated when the first author visited the Department
of Mathematics at the TU Darmstadt. He would like to thank for its
great hospitality.

\section{A generalized implicit function theorem for Banach spaces}

In this section we prove the generalized implicit function theorem.

\begin{theorem}\label{thm:implicit}
Let $X, Y$ and $Z$ be Banach spaces, $U\subeq X$, $V\subeq Y$ be
open $0$-neighborhoods. Let $\varphi:U\rightarrow V$ and
$\psi:V\rightarrow Z$ be $C^1$-maps.
Suppose\\
{\rm(1)} $\varphi(0)=0$, $\psi(0)=0$;\\
{\rm(2)} $\psi\circ\varphi\equiv0$;\\
{\rm(3)} $\im(d\varphi(0))=\ker(d\psi(0))$;\\
{\rm(4)} $\im(d\psi(0))$ is closed.\\
Then there exists a $0$-neighborhood $W\subeq V$ such that
$$\psi^{-1}(0)\cap W=\varphi(U)\cap W.$$
\end{theorem}

\begin{proof}
By the Open Mapping Theorem, there exist positive constants $C_1,
C_2$ such that for any $y\in\im(d\varphi(0))$ and
$z\in\im(d\psi(0))$, there are $\sigma(y)\in X$ and $\tau(z)\in Y$
with $d\psi(0)(\sigma(y))=y$, $d\psi(0)(\tau(z))=z$, and such that
$$\|\sigma(y)\|\leq C_1\|y\|,\qquad \|\tau(z)\|\leq C_2\|z\|.$$
Let $$C=\max\{C_1, C_2, \|d\varphi(0)\|, \|d\psi(0)\|, 1\}.$$ Choose
a $\delta>0$ with $B_X(0,\delta)\subeq U$, $B_Y(0,\delta)\subeq V$
such that
$$x\in
B_X(0,\delta)\Rightarrow\|d\varphi(x)-d\varphi(0)\|<\frac{1}{6C^3},$$
$$y\in
B_Y(0,\delta)\Rightarrow\|d\psi(y)-d\psi(0)\|<\frac{1}{6C^3},$$
where $B_X(0,\delta)$, $B_Y(0,\delta)$ are the open balls in $X$ and
$Y$ with centers $0$ and radius $\delta$. We claim that the
$0$-neighborhood $$W=B_Y\Big(0,\frac{\delta}{9C^3}\Big)$$ satisfies
the request of the theorem. In fact, we prove below that for $y\in
B_Y(0,\frac{\delta}{9C^3})$, if $\psi(y)=0$, then there exists $x\in
B_X(0,\frac{\delta}{2})$ such that $\varphi(x)=y$.

Suppose $y\in B_Y(0,\frac{\delta}{9C^3})$ such that $\psi(y)=0$.
Choose a point $x_0\in B_X(0,\frac{\delta}{162C^4})$, and define
$u_n\in Y$, $v_n, x_{n+1}\in X$ inductively by
$$\begin{cases}
u_n=\varphi(x_n)-y \ \ \ \text{(whenever $x_n\in U$)},\\
v_n=\sigma\big(u_n-\tau(d\psi(0)(u_n))\big),\\
x_{n+1}=x_n-v_n.\end{cases}$$
Here we use that $u_n - \tau(d\psi(0)(u_n)) \in \ker d\phi(0) = \im(d\psi(0))$.
We prove inductively that
$$\begin{cases}
\|x_n\|<\frac{\delta}{162}+\frac{20\delta}{81}\sum_{i=0}^{n-1}\frac{1}{2^i}
\ \ \ \text{(hence $x_n\in B_X(0,\frac{\delta}{2})\subeq U$)};\\
\|u_n\|<\frac{10\delta}{81C^3}\frac{1}{2^n};\\
\|v_n\|<\frac{20\delta}{81}\frac{1}{2^n} \leq \frac{10}{81} \delta < \frac{\delta}{2}.
\end{cases}$$

Since
$$\|x_0\|<\frac{\delta}{162C^4}\leq\frac{\delta}{162},$$
\begin{align*}
\|u_0\|&=\|\varphi(x_0)-y\|
\leq\|\varphi(x_0)-\varphi(0)\|+\|y\|\\
&<\|\int_0^1d\varphi(tx_0)(x_0)\ dt\|+\frac{\delta}{9C^3}
\leq\int_0^1\|d\varphi(tx_0)\|\|x_0\|\ dt+\frac{\delta}{9C^3}\\
&\leq\|x_0\|\int_0^1(\|d\varphi(0)\|+\|d\varphi(tx_0)-d\varphi(0)\|)\
dt+\frac{\delta}{9C^3}\\
&\leq\frac{\delta}{162C^4}\Big(C+\frac{1}{6C^3}\Big)+\frac{\delta}{9C^3}
<\frac{\delta}{81C^3}+\frac{\delta}{9C^3}
=\frac{10\delta}{81C^3},
\end{align*}
and
\begin{align*}
\|v_n\|&=\|\sigma(u_n-\tau(d\psi(0)(u_n)))\|\\
&\leq C\|u_n-\tau(d\psi(0)(u_n))\|
\leq C(\|u_n\|+\|\tau(d\psi(0)(u_n))\|)\\
&\leq C(\|u_n\|+C\|d\psi(0)(u_n)\|)
\leq C(\|u_n\|+C^2\|u_n\|)
\leq 2C^3\|u_n\|,
\end{align*}
which implies that $$\|v_0\|\leq2C^3\|u_0\|<\frac{20\delta}{81},$$
the inequalities hold for $n=0$. Suppose that the inequalities hold
for some $n\geq0$. We first have
$$\|x_{n+1}\|=\|x_n-v_n\|\leq\|x_n\|+\|v_n\|<\frac{\delta}{162}+\frac{20\delta}{81}\sum_{i=0}^{n}\frac{1}{2^i}.$$
Next we note that
\begin{align*}
u_{n+1}&=\varphi(x_{n+1})-y=\varphi(x_n-v_n)-y
=\varphi(x_n-v_n)-\varphi(x_n)+u_n\\
&=\varphi(x_n-v_n)-\varphi(x_n)+d\varphi(0)(v_n)+\tau(d\psi(0)(u_n)),
\end{align*}
and that $x_n, x_n - v_n \in B_X(0,\delta)$ implies that
$x_n - t v_n \in B_X(0,\delta)$ holds for $0 \leq t \leq 1$.
We thus obtain
\begin{align*}
\|u_{n+1}\|&\leq\|\varphi(x_n-v_n)-\varphi(x_n)+d\varphi(0)(v_n)\|+\|\tau(d\psi(0)(u_n))\|\\
&\leq\|\int_0^1(d\varphi(x_n-tv_n)-d\varphi(0))(v_n)\
dt\|+C\|d\psi(0)(u_n)\|\\
&\leq\int_0^1\|d\varphi(x_n-tv_n)-d\varphi(0)\|\|v_n\|\
dt\\
&\ \ \ \ +C\|\psi(\varphi(x_n))-\psi(y)-d\psi(0)(\varphi(x_n)-y)\|\\
&\leq\frac{1}{6C^3}\|v_n\|+C\|\int_0^1(d\psi(y+tu_n)-d\psi(0))(u_n)\ dt\|\\
&\leq\frac{1}{3}\|u_n\|+C\int_0^1\|d\psi(y+tu_n)-d\psi(0)\|\|u_n\|\
dt\\
&\leq\frac{1}{3}\|u_n\|+C\frac{1}{6C^3}\|u_n\|
\leq\frac{1}{2}\|u_n\|<\frac{10\delta}{81C^3}\frac{1}{2^{n+1}}.
\end{align*}
So we also have
$$\|v_{n+1}\|\leq2C^3\|u_{n+1}\|<\frac{20\delta}{81}\frac{1}{2^{n+1}}.$$
This proves the inequalities.

By the definition of $x_n$ and the inequalities, the sequence
$(x_n)_{n \in \N}$ is a Cauchy sequence, hence converges to some $x\in X$,
and we have
$$\|x\|=\|x_0-\sum_{i=0}^\infty
v_i\|\leq\|x_0\|+\sum_{i=0}^\infty\|v_i\|<
\frac{\delta}{162}+\frac{20\delta}{81}\sum_{i=0}^\infty\frac{1}{2^i}=\frac{\delta}{2}.$$
We also have
$\lim_{n\rightarrow\infty}u_n=0.$ Hence
$$y=\lim_{n\rightarrow\infty}(\varphi(x_n)-u_n)=\varphi(x).$$ This
proves the above claim, hence the theorem.
\end{proof}

\begin{remark}
Theorem~\ref{thm:implicit} does not hold without the assumption that
$\im(d\psi(0))$ being closed. Here is a counter-example. Let $H$ be
a Hilbert space with an orthonormal basis $\{e_n:n\in\N\}$, $f$ be
the $C^1$-map from $H$ to itself defined by $f(v)=\|v\|v$, $A$ be
the linear map from $H$ to itself determined by
$Ae_n=-\frac{1}{n}e_n$. Let $X=0$, $Y=Z=H$, $\phi\equiv0$,
$\psi=f+A$. Then $d\psi(0)=A$, which is injective. It is obvious
that all the conditions in Theorem~\ref{thm:implicit} hold except
that $\im(d\psi(0))=\im(A)$ is closed.
$\frac{1}{n}e_n\in\psi^{-1}(0)$, but $\frac{1}{n}e_n$ does not
belong to $\im(\varphi)=\{0\}$ for each $n\in\N$.
\end{remark}

\begin{corollary}
Let $Y$ and $Z$ be Banach spaces, $V\subeq Y$ an open
$0$-neigh\-bor\-hood, $\psi:V\rightarrow Z$ be a $C^1$-map and $X :=
\ker(d\psi(0))$. Suppose $\psi(X \cap V)=\{0\}$ and that
$\im(d\psi(0))$ is closed. Then there exists a $0$-neighborhood
$W\subeq V$  such that
$$\psi^{-1}(0)\cap W=X\cap W.$$
\end{corollary}

\section{Parametrization of fibers of differentiable maps}

In this section we refine the Implicit Function
Theorem~\ref{thm:implicit} in such a way that it also provides a manifold
structure on the $0$-fiber of $\psi$. Clearly, this is not a submanifold
in the sense of Bourbarki~\cite{Bou89} because its tangent space
need not be complemented, but it nevertheless is a Banach manifold.

\begin{proposition} \label{prop:open}
Let $X$, $Y$ and $Z$ be Banach spaces. Then the following
assertions hold:
\begin{description}
\item[\rm(1)] The set $L_s(X,Y)$ of all surjective linear operators
from $X$ to $Y$ is open.
If $f \in L(X,Y)$ and $B_Y(0,r) \subeq f(B_X(0,1))$, then
$f+g \in L_s(X,Y)$ whenever $\|g\| < r$.
\item[\rm(2)] The set $L_e(X,Y)$ of all closed embeddings of
$X$ into $Y$ is open.
If $f \in L(X,Y)$ satisfies
$\|f(x)\| \geq r \|x\|$ for some $r > 0$, then
$f+g \in L_e(X,Y)$ whenever $\|g\| < r$.
\item[\rm(3)] The set of pairs $(f,g)$ with $\im(f) = \ker(g)$
is an open subset in
$$ \{ (f,g) \in L_e(X,Y) \times L_s(Y,Z) \: gf = 0\}. $$
\end{description}
\end{proposition}

\begin{proof} (1) (This proof is due to H.~Gl\"ockner. Another
argument is given in \cite{Bro76}, Lemma and Theorem 2.1 on pp.20/21.)
Let $f \: X \to Y$ be open. In view of the Open Mapping Theorem,
there exists an $r > 0$ with $B_Y(0,r) \subeq f(B_X(0,1))$.
We claim that $f + g$ is surjective for each $g \in L(X,Y)$ with
$\|g\| < r$.

Let $v \in Y$ and $a := \frac{\|g\|}{r} < 1$. Then there exists
some $w_0 \in X$ with $f(w_0) = v$ and $\|w_0\| \leq \frac{\|v\|}{r}$.
Then we have
$(f+g)(w_0) = v + g(w_0)$ with
$$ \|g(w_0)\| \leq r^{-1}\|v\|\|g\| = a\|v\|. $$
Now we pick $w_1 \in X$ with $f(w_1) = - g(w_0)$
and $\|w_1\| \leq r^{-1}\|g(w_0)\|\leq a\|w_0\|\leq \frac{a}{r}\|v\|$
and obtain
$$ (f+g)(w_0 + w_1) = v + g(w_1)
\quad \mbox{ with } \quad
\|g(w_1)\| \leq a^2 \|v\|. $$
Iterating this procedure, we obtain a sequence $(w_n)_{n \in \N}$
in $X$ with
$\|w_n\| \leq a^n\frac{\|v\|}{r}$, so that
$w := \sum_{i = 0}^\infty w_i$ converges, and
$f(w) = \lim_{n \to \infty} v + g(w_n) = v$.

(2) The Open Mapping Theorem implies that $f \in L(X,Y)$ is a closed embedding
if and only if there exists an $r > 0$ with
$\|f(x)\| \geq r\|x\|$ for all $x \in X$.
If this is the case and $g \in L(X,Y)$ satisfies $\|g\| < r$, then
$$ \|(f+g)(x)\| \geq \|f(x)\|-\|g(x)\| \geq (r-\|g\|)\|x\| $$
implies that $f + g$ is a closed embedding.

(3) Let $(f,g) \in L_e(X,Y) \times L_s(Y,Z)$. Without loss of
generality, we may assume that $f$ is an isometric embedding and
that $g \: Y \to Z$ is a metric quotient map, i.e., $\|g(x)\| =
\inf_{y \in x + \ker g} \|y\|$. We show that if $\delta <
\frac{1}{10}$ then $\|f-\tilde f\|, \|g-\tilde g\| < \delta$ and
$\tilde g\tilde f =0$ imply that $\im(\tilde f) = \ker(\tilde g)$.

Let $v \in \ker\tilde g$. We put $v_1 := v$ and assume that we have
already constructed a sequence $v_1,\ldots, v_n$ in $\ker \tilde g$
with
$$ \|v_{j+1}\| \leq \frac{1}{2} \|v_j\| \quad \mbox{ and } \quad
v_{j+1} - v_j \in \im(\tilde f), \quad j =1,\ldots, n-1. $$

To find $v_{n+1}$, we first observe that
$$ \|g(v_n)\| = \|g(v_n) - \tilde g(v_n)\|
\leq \|g-\tilde g\| \|v_n\| \leq \delta \|v_n\|. $$
Hence there exists $w_n \in Y$ with $\|w_n\| \leq 2 \delta \|v_n\| \leq \|v_n\|$
and $g(w_n) = g(v_n)$. Then $v_n - w_n \in \ker g = \im f$ and
$x_n := f^{-1}(v_n - w_n)$ satisfies
$$\|x_n\| = \|v_n - w_n\| \leq 2 \|v_n\|. $$
For $v_{n+1} := v_n - \tilde f(x_n)$ we now have
$v_{n+1} \in \ker \tilde g$ and
\begin{eqnarray*}
\|v_{n+1}\|
&=& \|v_n - \tilde f(x_n)\| = \|v_n  - f(x_n) + (f - \tilde f)(x_n)\| \\
&=& \|w_n  + (f - \tilde f)(x_n)\| \leq \|w_n\| + \delta \|x_n\| \\
&\leq& 2\delta \|v_n\| + 2 \delta \|v_n\| \leq \frac{1}{2} \|v_n\|.
\end{eqnarray*}
We thus obtain sequences $(v_n)$ and $(x_n)$ satisfying the above
conditions. In particular, $x := \sum_{n=1}^\infty x_n$ converges in
$F$ and $\tilde f(x) = \sum_{n=1}^\infty v_{n} - v_{n+1} = v_1 = v$.
\end{proof}

With the preceding proposition we can strengthen Theorem
\ref{thm:implicit} under the assumption that $d\varphi(0)$ is
injective and $d\psi(0)$ is surjective as follows.

\begin{proposition} \label{rem-open}
Let $X, Y$ and $Z$ be Banach spaces, $U\subeq X$, $V\subeq Y$ be
open $0$-neighborhoods. Let $\varphi:U\rightarrow V$ and
$\psi:V\rightarrow Z$ be $C^1$-maps.
Suppose\\
{\rm(1)} $\varphi(0)=0$, $\psi(0)=0$;\\
{\rm(2)} $\psi\circ\varphi\equiv0$;\\
{\rm(3)} $\im(d\varphi(0))=\ker(d\psi(0))$;\\
{\rm(4)} $d\varphi(0)$ is injective, $d\psi(0)$ is surjective.\\
Then there exists an open $0$-neighborhood $U' \subeq U$ such that
$\phi\res_{U'}$ is an open map into $\psi^{-1}(0)$.
\end{proposition}

\begin{proof}
For each $u \in U$ the relation $\psi \circ \phi = 0$ implies that
$$ d\psi(\phi(u)) d\phi(u)= 0, $$
so that Proposition~\ref{prop:open} implies that
$$ \ker(d\psi(\phi(u))) = \im(d\phi(u)) $$
holds for each $u$ in some open $0$-neighborhood $U' \subeq U$.
Applying Theorem~\ref{thm:implicit} to the maps $\psi_u(y) :=
\psi(\phi(u)+y)$ and $\phi_u(x) := \phi(u+x)$, it follows that
$\phi\res_{U'}$ is an open map to $\psi^{-1}(0)$.
\end{proof}

\begin{lemma} \label{lem:quasiiso} Let $U \subeq X$ be an open subset containing $0$,
$f \: U \to Y$ a $C^1$-map for which $df(0)$ is a closed embedding.
Then there exists an open $0$-neighborhood $U'\subeq U$ and
constants $0 < c < C$ with
$$ c\|x- y\| \leq \|f(x)-f(y)\| \leq C\|x-y\| \quad \mbox{ for } \quad x,y \in
U'. $$
\end{lemma}

\begin{proof} Since $df(0)$ is a closed embedding, there exists
some $r > 0$ with $\|df(0)x\| \geq r\|x\|$ for $x \in X$. Pick an
open ball $U'\subeq U$ around $0$ with
$$\|df(p)-df(0)\| \leq \frac{r}{4} \quad \mbox{ for } \quad p \in U'$$
and note that
this implies that
$$ \|df(p)v\| \geq \frac{3r}{4}\|v\| \quad \mbox{ for } \quad v \in X. $$
We now obtain for $x,y \in U'$:
\begin{eqnarray*}
&& \ \ \ \|f(y)-f(x) - df(x)(y-x)\| \\
&=&\| \int_0^1 df(x + t(y-x))(y-x)\, dt - df(x)(y-x)\|\\
&=&\| \int_0^1 \big(df(x + t(y-x))-df(x)\big)(y-x)\, dt\|
 \leq \frac{r}{2}\|y-x\|,
\end{eqnarray*}
which in turn leads to
$$ \|f(y)-f(x)\|
\geq \|df(x)(y-x)\| - \frac{r}{2}\|y-x\| \geq \frac{3r}{4}\|y-x\| -
\frac{r}{2}\|y-x\|  = \frac{r}{4}\|y-x\|. $$ With $c := \frac{r}{4}$
and $C := \|df(0)\| + \frac{r}{4}$ the assertion now follows from
$$ \|f(y)-f(x)\|
= \| \int_0^1 df(x + t(y-x))(y-x)\, dt\|
\leq C \|y-x\|. $$
\end{proof}

Let $\psi \: V \to Z$ be a $C^k$-map from an open $0$-neighborhood
$V$ in a Banach space $Y$ into another Banach space $Z$ with
$\psi(0)=0$ and $d\psi(0)$ is surjective. Assume that there is an
open $0$-neighborhood $U$ in some Banach space $X$ and a $C^k$-map
$\phi \: U \to \psi^{-1}(0) \subeq V$ such that $\phi(0) = 0$,
$d\phi(0)$ is injective, and $\im(d\phi(0)) = \ker(d\psi(0)$. Then
by Proposition~\ref{rem-open} and Lemma~\ref{lem:quasiiso},
$\psi^{-1}(0)$ a locally a Banach $C^k$-manifold around $0$, and a
local $C^k$-chart is provided by $\phi$. To get the global version
of such result, we need more preparation.

\begin{lemma} \label{lem:cont} The map
$$ \mu \: \{ (f,g) \in L_e(X,Y)^2 \: \im(f) = \im(g) \} \to \GL(X), \quad
(f,g) \mapsto f^{-1} \circ g $$
is continuous.
\end{lemma}

\begin{proof} Fix $(f_0, g_0)\in L_e(X,Y)^2$ with
$\im(f_0) = \im(g_0)$. We may w.l.o.g.\ assume that $f_0$ is an isometric
embedding $X \into Y$. Then for each $h \in L(X,Y)$ with
$\|h\| < 1$ the map $f_0 + h$ is a closed embedding with
$$ \|(f_0+h)v\|\geq (1 - \|h\|) \|v\| \quad \mbox{ for } \quad
v \in X. $$
We conclude that for any $g \in L(X,Y)$ with $\im(g) \subeq \im(f_0+h)$ we have
$$ \|(f_0+h)^{-1}g\| \leq (1 - \|h\|)^{-1}\|g\|. $$
For any $h$, $g \in L(X,Y)$ with
$\im(f_0 + h) = \im(g)$ and $x \in X$ we thus obtain
\begin{eqnarray*}
\|(f_0+h)^{-1}g- f_0^{-1}g_0\|
&=& \| (f_0+h)^{-1}\big(g - (f_0 + h)f_0^{-1}g_0\big)\| \\
&\leq& (1-\|h\|)^{-1} \|g - (f_0 + h)f_0^{-1}g_0\| \\
&\leq &(1-\|h\|)^{-1} \big(\|g - g_0\| + \|\1 - (f_0 + h)f_0^{-1}\|
\|g_0\|\big).
\end{eqnarray*}
Since $f_0$ is assumed to isometric, the map
$$\1 - (f_0 + h)f_0^{-1} \: \im(f_0) = \im(g_0) \to Y $$
satisfies
$$ \|\1 - (f_0 + h)f_0^{-1}\| = \|f_0 - (f_0 + h)\| = \|h\|. $$
Therefore the estimate above implies the continuity of
$\mu$.
\end{proof}

\begin{theorem}\label{thm:ck} Let $X, Y$, and $Z$ be Banach spaces.
Let $V \subeq Y$ be an open subset containing $0$ and $\psi \: V \to
Z$ be a $C^k$-map, $k \in \N \cup \{\infty\}$, such that $\psi(0) =
0$ and $d\psi(0)$ is surjective. Further, let $U_1, U_2 \subeq X$ be
open subsets containing $0$ and $\phi_i \: U_i \to S :=
\psi^{-1}(0)$ be $C^k$-maps such that $\phi_i(0) = 0$, $d\phi_i(0)$
is injective, and $\im(d\phi_i(0)) = \ker d\psi(0)$ for $i = 1,2$.
Then there exist open $0$-neighborhoods $U_i' \subeq U_i$ such that
$\phi_i\res_{U_i'}$ are homeomorphisms onto a $0$-neighborhood in
$S$ and the map
$$ \phi_{12} := \phi_1^{-1} \circ \phi_2 \: U_2' \cap \phi_2^{-1}(\phi_1(U_1')) \to X $$
is $C^k$.
\end{theorem}

\begin{proof}
We prove the theorem by induction over $k \in \N$. If it holds for
each $k \in \N$, it clearly holds for $k = \infty$.

We first prove the theorem for $k=1$. Using
Proposition~\ref{rem-open}, we may w.l.o.g.\ assume that $\phi_i$ is
a homeomorphism onto an open $0$-neighborhood in $S$ and that
$d\phi_i(x)$ is injective, $d\psi(\phi_i(x))$ is surjective, and
$\im(d\phi_i(x)) = \ker(d\psi(\phi_i(x)))$ hold for all $x \in U_i$.
Next we use Lemma~\ref{lem:quasiiso} to see that we may further
assume that there exist constants $0 < c < C$ such that for $x,y \in
U_i$ we have
\begin{eqnarray} \label{eq:esti}
c\|x- y\| \leq \|\phi_i(x)-\phi_i(y)\| \leq C\|x-y\|.
\end{eqnarray}

Now fix $p \in U_2'$ with $\phi_2(p) \in \phi_1(U_1)$.
We have to show that $\phi_{12}$ is differentiable
in $p$. Let $y := \phi_2(p)$ and $q := \phi_1^{-1}(y)$.
Since $d\phi_2(p)$ and $d\phi_1(q)$ are closed embeddings, the
linear map $A := d\phi_1(q)^{-1} \circ d\phi_2(p) \: X \to X$ is
invertible.

We have
$$ \phi_2(p+h)
= \phi_2(p)+ d\phi_2(p)h + r_2(h)
=  y + d\phi_2(p)h + r_2(h) $$
with $\lim_{h \to 0} \frac{r_2(h)}{\|h\|} = 0$ and
$$ \phi_1(q+h)
= \phi_1(q)+ d\phi_1(q)h + r_1(h) = y+ d\phi_1(q)h + r_1(h) $$
with $\lim_{h \to 0} \frac{r_1(h)}{\|h\|} = 0$.
For $\nu(h) := \phi_{12}(p+h) - q$ we then have
$$  y + d\phi_2(p)h + r_2(h)
= \phi_2(p+h) = \phi_1(q+\nu(h))
= y+ d\phi_1(q)\nu(h) + r_1(\nu(h)), $$
so that
$$  d\phi_2(p)h + r_2(h) =  d\phi_1(q)\nu(h) + r_1(\nu(h)). $$
This implies that $r_2(h) - r_1(\nu(h)) \in \ker d\psi(y) = \im(d\phi_1(q))$
and that
$$  \nu(h) = Ah + (d\phi_1(q))^{-1}(r_2(h)-r_1(\nu(h))). $$
From (\ref{eq:esti}) we derive
$$ \frac{c}{C} \|h\| \leq \|\nu(h)\| \leq  \frac{C}{c} \|h\|, $$
and hence
$$ \lim_{h \to 0} \frac{r_1(\nu(h))}{\|h\|} = 0. $$
Therefore $\nu$ is differentiable in $0$ with
$d\nu(0)=A$.

We conclude that $\phi_{12}$ is differentiable in
$p$ with
$$ d\big(\phi_{12}\big)(p)
=  (d\phi_1(q))^{-1} \circ d\phi_2(p). $$ The continuity of the
differential $d(\phi_{12})$ follows from Lemma~\ref{lem:cont}. This
finishes the proof of the theorem for $k=1$.

Now let us assume that $k > 1$ and that the theorem holds for
$C^{k-1}$-maps, showing that $\phi_{12}$ is a $C^{k-1}$-map. It
remains to show that the tangent map $T\phi_{12}$ is also $C^{k-1}$,
which therefore proves that $\phi_{12}$ is $C^k$.

The tangent map
$$ T\psi \: TV \cong V \times Y \to TZ \cong Z \times Z, \quad
(x,v) \mapsto (\psi(x),d\psi(x)v)$$
is a $C^{k-1}$-map
with $T\psi(0,0) = (0,0)$
for which
$d(T\psi)(0,0) \cong d\psi(0) \times d\psi(0)$ is surjective.
Further, the maps
$$T\phi_i \: TU_i \cong U_i \times X \to
\tilde S := (T\psi)^{-1}(0,0)$$
are $C^{k-1}$-maps with
$$ T\phi_i(0,0) = (0,0) \quad \mbox{ and } \quad
\im(d(T\phi_i)(0,0)) = \ker d(T\psi)(0,0), i = 1,2. $$
Our induction hypothesis implies that $T\phi_{12} = T\phi_1^{-1} \circ T\phi_2$
is a $C^{k-1}$-map, and this completes the proof.
\end{proof}

\begin{theorem}\label{thm:regular} {\rm(Regular Value Theorem; without complements)}
Let $\psi \: M \to N$ be a $C^k$-map bet\-ween Banach
$C^k$-manifolds, $q \in N$, $S := \psi^{-1}(q)$ and $X$ be a Banach
space. Assume that
\begin{description}
\item[\rm(1)]  $q$ is a regular value in the sense that for each
$p \in S$ the differential $d\psi(p) \: T_p(M) \to T_q(N)$ is
surjective, and
\item[\rm(2)] for every $p\in S$, there exists an open $0$-neighborhood
$U_p \subeq X$ and a $C^k$-map $\phi_p \: U_p \to S \subeq M$ with
$\phi_p(0) = p$, $d\phi_p(0)$ is injective, and $\im(d\phi_p(0)) =
\ker(d\psi(p)$.
\end{description}
Then $S$ carries the structure of a Banach $C^k$-manifold
modeled on $X$.
\end{theorem}

\begin{proof} In view of Proposition~\ref{rem-open},
 the maps $\phi_p$
whose existence is assured by (2) yield local charts and
Theorem~\ref{thm:ck} implies that the chart changes are also $C^k$,
so that we obtain a $C^k$-atlas of $S$.
\end{proof}

\section{Applications to Lie subgroups}

If $G_1$ and $G_2$ are Banach--Lie groups, then the existence of
canonical coordinates (given by the exponential function) implies
that each continuous homomorphism $\phi \: G_1 \to G_2$ is actually
smooth, hence a morphism of Lie groups. This implies in particular
that a topological group $G$ carries at most one Banach--Lie group
structure.

If $G$ is a Banach--Lie group, then we call a closed subgroup $H \subeq G$
a {\it Banach--Lie subgroup} if it is a Banach--Lie group with respect to the
subspace topology. It is an important problem in infinite-dimensional
Lie theory to find good criteria for a closed subgroup of a Banach--Lie
group to be a Banach--Lie subgroup.

Regardless of any Lie group structure, we may defined for
each closed subgroup $H$ of a Banach--Lie group its
Lie algebra by
$$ \L(H) := \{ x \in \L(G) \: \exp_G(\R x)\subeq H\}, $$
and the Product and Commutator Formula easily imply that
$\L(H)$ is a closed Lie subalgebra of the Banach--Lie algebra $\L(G)$
(cf.\ \cite{Ne06}, Lemma~IV.3.1). By restriction we thus obtain an
exponential map
$$ \exp_H \: \L(H) \to H. $$
From \cite{Ne06}, Thm.~IV.4.16 we know that there exists a Banach--Lie group
$H_L$ with Lie algebra $\L(H)$ and an injective morphism of Lie groups
$i_H \: H_L \into G$ for which $i_H(H_L) = \la \exp_G \L(H) \ra$ is the
subgroup of $H$ generated by the image of $\exp_H$.
Since $\L(H)$ is a closed subspace of $\L(G)$, the map
$\L(i_H) \: \L(H) \to \L(G)$ is a closed embedding, and
$H$ is a Banach--Lie subgroup if and only if $i_H \: H_L \to H$ is an open map,
which is equivalent to the existence of an open $0$-neighborhood
$U \subeq \L(H)$ for which $\exp_H\res_U \: U \to H$ is an open embedding
(cf.\ \cite{Ne06}, Thm.~IV.3.3).
This notion of a Banach--Lie subgroup
is weaker than the concept used by Bourbaki (\cite{Bou89}), where it is
assumed that $\L(H)$ has a closed complement, a requirement that is
not always satisfied and hard to verify in many concrete situations.

For closed normal subgroups $H$ of Banach--Lie groups one has the nice
characterization that $H$ is a Banach--Lie subgroup if and only if the topological
quotient group $G/H$ is a Banach--Lie group (\cite{GN03}). It would be nice
to have a similar criterion for general subgroups:

\begin{conjecture} \label{conj}
A closed subgroup $H$ of a Banach--Lie group $G$ is a
Banach--Lie subgroup if and only if $G/H$ carries the structure of a
Banach manifold for which the quotient map
$q \: G \to G/H, g \mapsto gH$ has surjective differential in each point
and $G$ acts smoothly on $G/H$.
\end{conjecture}

The main difficulty arises from the possible absence of closed complements
to the closed subspace $\L(H)$ in $\L(G)$. If $H$ is a Banach--Lie subgroup of $G$
and $\L(H)$ has a closed complement,
then the classical Inverse Function Theorem provides natural charts on the
quotient space $G/H$, hence a natural manifold structure with all nice
properties (cf.\ \cite{Bou89}, Ch.~3). Without a closed complement for $\L(H)$
it is not known how to construct charts of $G/H$. However,
the natural model space is the quotient space $\L(G)/\L(H)$.

With the aid of the ``implicit function theorem'', we can now prove
one half of the conjecture above:

\begin{theorem} \label{thm:liesub} Let $G$ be a Banach--Lie group and
$\sigma \: G \times M \to M$ a smooth action of $G$ on a Banach
manifold $M$, $p\in M$. Suppose that for the derived action
$$ \dot\sigma \: \L(G) \to {\cal V}(M), \quad
\dot\sigma(x)(m) := -d\sigma(\1,m)(x,0) $$ the subspace
$\dot\sigma(\L(G))(p)$ of $T_p(M)$ is closed. Then the stabilizer
$$ G_p := \{ g \in G \: g.p = p\} $$
is a Lie subgroup.
\end{theorem}

\begin{proof} First we note that $G_p$ is a closed subgroup of
$G$ with
$$ \L(G_p) = \{ x \in \L(G) \: \dot\sigma(x)(p) =0\}. $$
We consider the smooth maps
$$ \psi \: G \to M, \quad g \mapsto g.p
\quad \mbox{ and } \quad
\phi \: \L(G_p) \to G_p, \quad x \mapsto \exp_G(x). $$
Then $G_p = \psi^{-1}(p)$, $\phi(0)= \1$,
$$ \im(d\phi(0)) = \L(G_p) = \ker(d\psi(\1)). $$
Since $\im(d\phi(0))=\dot\sigma(\L(G))(p)$ is closed, Theorem
\ref{thm:implicit} implies that there exists arbitrarily small open
$0$-neighborhood $U \subeq \L(G_p)$ such that $\exp_G(U)$ is a
$\1$-neighborhood in $G_p$. By \cite{Ne06}, Thm.~IV.3.3 (cf.\ also
\cite{Hof75},  Prop.~3.4), $G_p$ is a Banach--Lie subgroup of $G$.
\end{proof}

\begin{corollary} Let $G$ be a Banach--Lie group and
$H \subeq G$ a closed subgroup for which $G/H$ carries the structure
of a Banach manifold. Suppose the quotient map $q \: G \to G/H, g
\mapsto gH$ has surjective differential in some point $g_0\in G$ and
$G$ acts smoothly on $G/H$. Then $H$ is a Banach--Lie subgroup of~$G$.
\end{corollary}

\begin{proof} The action of $G$ on $G/H$ is given by
$\sigma(g,xH) = gxH = q(gx),$
so that
$$ \dot\sigma(\L(G))(gH) = dq(g)(\L(G)) = T_{gH}(G/H) $$
and the stabilizer of $p := gH = q(g) \in G/H$ is $gHg^{-1}$. By
Theorem~\ref{thm:liesub}, $g_0 Hg_0^{-1}$ is a Banach--Lie subgroup. So
$H$ is a Banach--Lie subgroup.
\end{proof}

\begin{remark} \label{rem:more-gen}
The preceding theorem has a natural generalization to the following
setting. Let $H$ be a closed subgroup of a Banach--Lie group $G$.
Suppose that there exists an open $\1$-neighborhood $U_G \subeq G$
and a $C^1$-map $F \: U_G \to M$ to some Banach manifold $M$ such
that $dF(\1)$ is surjective and that for $m := F(\1)$ we have
$$F^{-1}(m) \cap U_G = H \cap U_G. $$
Then $H$ is a Banach--Lie subgroup of $G$.
\end{remark}

\begin{problem} Suppose that the Banach--Lie group $G$ acts smoothly
on the Banach manifold $M$ and $p \in M$. Is it always true that the
stabilizer $G_p$ is a Lie subgroup of $G$?

From the theory of algebraic Banach--Lie groups it follows that
for each Banach space $Y$, for each $y \in Y$ and each closed subspace
$X \subeq Y$, the subgroups
$$ \GL(Y)_y := \{ g \in \GL(Y) \: g(y)=y\} $$
and
$$ \GL(Y)_X := \{ g \in \GL(Y) \: g(X)=X\} $$
are Banach--Lie subgroups. Since inverse images of Lie subgroups are
Lie subgroups (\cite{Ne04}, Lemma IV.11)  the problem
has a positive solution for linear actions and for stabilizers of
closed subspaces.

If $X$ is not complemented in $Y$, then it is not clear how to turn
the Grassmannian $\Gr_X(Y) := \{ g.X \: g \in \GL(Y) \}$ into a smooth manifold
on which $\GL(Y)$ acts. We are therefore in a situation where we know that
the closed subgroup $\GL(X)_Y$ is a Banach--Lie subgroup, but nothing is known
on the quotient $\Gr_X(Y) \cong \GL(Y)/\GL(Y)_X$. Therefore a solution of
Conjecture~\ref{conj}
would in particular lead to a natural Banach manifold structure
on the Grassmannian of all closed subspaces of $Y$.
\end{problem}

\section{Weil's local rigidity for Banach--Lie groups}

Let $G$ be a topological group, $\Gamma$ a finitely generated (discrete)
group.
Let $\Hom(\Gamma,G)$ be the set of all homomorphisms from $\Gamma$
to $G$, endowed with the compact-open topology, which coincides with the
topology of pointwise convergence. A homomorphism
$r\in\Hom(\Gamma,G)$ is \emph{locally rigid} if there is an open
neighborhood $U$ of $r$ in $\Hom(\Gamma,G)$ such that for any $r'\in
U$, there exists some $g\in G$ with
$r'(\gamma)=gr(\gamma)g^{-1}$ for all $\gamma\in\Gamma$.

If, moreover, $G$ is a Banach--Lie group, then every homomorphism
$r\in\Hom(\Gamma,G)$ defines naturally a $\Gamma$-module structure
on the Lie algebra $\L(G)$ of $G$, via $\gamma.v=\Ad(r(\gamma))(v)$,
$\gamma\in\Gamma, v\in\L(G)$. For a fixed $r\in\Hom(\Gamma,G)$
and a map $f \: \Gamma \to G$ the pointwise product $f\cdot r \: \Gamma \to G$
is a homomorphism if and only if $f$ satisfies the cocycle condition
$$ f(\gamma_1\gamma_2) = f(\gamma_1) \cdot (r(\gamma_1)f(\gamma_2)r(\gamma_1)^{-1}). $$
For a smooth family of homomorphisms $r_t\in\Hom(\Gamma,G)$
$(-\varepsilon<t<\varepsilon)$ with $r_0=r$, this observation
directly implies that
the differential of $r_tr^{-1}$ at $t=0$ is a cocycle of $\Gamma$ in
$\L(G)$, and if $r_t=g_trg_t^{-1}$ for a smooth curve $g_t$ in $G$
with $g_0=e$, then the above differential is a coboundary. So it is
reasonable to expect some relation between the local rigidity of $r$
and the vanishing of $H^1(\Gamma,\L(G))$.

The following classical result is due to Weil \cite{We64}.

\begin{theorem}\label{thm:Weil}{\rm(Weil)}
Let $G$ be a finite-dimensional Lie group and $r\in$ \break $\Hom(\Gamma,G)$. If
$H^1(\Gamma,\L(G))=\{0\}$, then $r$ is locally rigid.
\end{theorem}

The main ingredient in Weil's proof of Theorem~\ref{thm:Weil} is the
finite-dimensional version of Theorem~\ref{thm:implicit}. Using
Theorem~\ref{thm:implicit}, in this section we generalize Weil's
Theorem to Banach--Lie groups. To state our result, we first recall
a definition of the first cohomology of groups, which is related to
presentations of groups.

Let $\Gamma$ be a group. Let $$1\rightarrow
R\stackrel{\iota}\rightarrow F\stackrel{\pi}\rightarrow
\Gamma\rightarrow1$$ be a presentation of $\Gamma$, where $F$ is a
free group on a subset $S\subeq F$, $R$ is a free group on a subset
$T\subeq R$. For $t\in T$, $\iota(t)$ can be expressed as a reduced
word
$$\iota(t)=s_{t,1}^{\epsilon_{t,1}}\cdots s_{t,m_t}^{\epsilon_{t,m_t}}$$ in
$S$, where $s_{t,j}\in S$, $\epsilon_{t,j}=\pm1$.

Let $V$ be a $\Gamma$-module. For a set $A$, let $\Map(A,V)$ denote
the set of all maps from $A$ to $V$. Then $\Map(A,V)$ has a natural
abelian group structure induced from that of $V$. We define the
coboundary operators $$\delta^0:V\rightarrow\Map(S,V) \quad \mbox{ and } \quad
\delta^1:\Map(S,V)\rightarrow\Map(T,V)$$ as follows:
$$\delta^0(v)(s)=v-\pi(s).v,\quad \delta^1(c)(t)=\sum_{j=1}^{m_t}\epsilon_{t,j}\pi(s_{t,1}^{\epsilon_{t,1}}\cdots
s_{t,j-1}^{\epsilon_{t,j-1}}s_{t,j}^{\epsilon'_{t,j}}).c(s_{t,j})$$
for $v\in V$ and $c\in\Map(S,V)$, where
$$\epsilon'_{t,j}=\begin{cases}
0, &\epsilon_{t,j}=1;\\
-1, &\epsilon_{t,j}=-1.
\end{cases}$$
Then $\delta^0$ and $\delta^1$ are homomorphisms of abelian groups,
$\delta^1\circ\delta^0=0$. The first cohomology group of $\Gamma$ in
$V$ is, by definition, $H^1(\Gamma,V)=\ker(\delta^1)/\im(\delta^0)$.

\begin{remark}\label{r:coincide}
This definition of $H^1(\Gamma,V)$ coincides with the usual one
which is defined by $H^1(\Gamma,V)=Z^1(\Gamma,V)/B^1(\Gamma,V)$,
where
$$Z^1(\Gamma,V)=\{f\in\Map(\Gamma,V)|(\forall \gamma_1, \gamma_2 \in \Gamma)\
f(\gamma_1\gamma_2)=f(\gamma_1)+\gamma_1.f(\gamma_2)\},$$
$$B^1(\Gamma,V)=\{f\in\Map(\Gamma,V)|(\exists v \in V)(\forall \gamma\in\Gamma)\
f(\gamma)=v-\gamma.v \}.$$
A direct way to see this is to look at the bijection between
$Z^1(\Gamma,V)$ and $\ker(\delta^1)$ which sends $f\in Z^1(\Gamma,V)$
to the element $s\mapsto f(\pi(s))$ in $\ker(\delta^1)$, under which
$B^1(\Gamma,V)$ is mapped bijectively onto $\im(\delta^0)$. A more
intrinsic way to see the coincidence is explained in Remark
\ref{r:hausdorff} below.
\end{remark}

Now suppose that $\Gamma$ is a finitely generated group. Then the
presentation of $\Gamma$ may be chosen in such a way that $F$ is finitely
generated. So we may assume that $S$ is finite, $T$ is countable.
Suppose moreover that $V$ is a continuous Banach $\Gamma$-module.
Then for a countable set $A$, $\Map(A,V)$ a Fr\'{e}chet space with
respect to the semi-norms $$p_a(f)=\|f(a)\|,\quad a\in A.$$ If,
moreover, $A$ is finite, then $\Map(A,V)$ is a Banach space with
respect to the norm $$\|f\|_A=\sup_{a\in A}\|f(a)\|.$$ In
particular, $(\Map(S,V), \|\cdot\|_S)$ is a Banach space,
$(\Map(T,V),\{p_t(\cdot)|t\in T\})$ is a Fr\'{e}chet space. It is
easy to see that the coboundary operators $\delta^0$ and $\delta^1$
are continuous. We say that $\delta^1$ \emph{has closed image} if
$\im(\delta^1)$ is a closed subspace of $\Map(T,V)$. Our
generalization of Weil's Theorem can be stated as follows.

\begin{theorem}\label{thm:rigidity}
Let $G$ be a Banach--Lie group, $\Gamma$ be a finitely generated
group and $r\in\Hom(\Gamma,G)$. If $H^1(\Gamma,\L(G))=\{0\}$ and $\delta^1$
has closed image, then $r$ is locally rigid.
\end{theorem}

In general, the set of generators $T$ of the group $R$ is infinite.
So $\Map(T,V)$ is only a Fr\'{e}chet space. The following lemma
enable us to reduce the problem to the category of Banach spaces so
as to we can apply Theorem~\ref{thm:implicit}.

\begin{lemma}\label{lem:finite}
Let $\Gamma$ be a finitely generated group, $V$ a continuous Banach
$\Gamma$-module, and let the notation be as above. If $\delta^1$ has
closed image, then there exists a finite subset $T'\subeq T$ such
that the map $\delta^1_{T'}:\Map(S,V)\rightarrow\Map(T',V)$ defined
by $\delta^1_{T'}(c)=\delta^1(c)|_{T'}$ satisfies
$\ker(\delta^1_{T'})=\ker(\delta^1)$, and $\im(\delta^1_{T'})$ is
closed in $(\Map(T',V),\|\cdot\|_{T'})$.
\end{lemma}

\begin{proof}
Since $\delta^1$ is continuous, $\ker(\delta^1)$ is a closed
subspace of $\Map(S,V)$. So $X=\Map(S,V)/\ker(\delta^1)$ is a Banach
space with respect to the norm $$\|[c]\|'=\inf_{z\in
\ker(\delta^1)}\|c+z\|_S,$$ where $[c]=c+\ker(\delta^1)$, and the
linear injection $\widetilde{\delta^1}:X\rightarrow\Map(T,V)$
induced by $\delta^1$ is continuous. Since
$\im(\widetilde{\delta^1})=\im(\delta^1)$ is closed in $\Map(S,V)$,
$\im(\widetilde{\delta^1})$ is a Fr\'{e}chet space with respect to
the semi-norms $p_t$, $t\in T$. By the Open Mapping Theorem for
Fr\'{e}chet spaces (see, for example, \cite{Bo87}, Section I.3,
Corollary 1), $\widetilde{\delta^1}$ is an isomorphism of
Fr\'{e}chet spaces from $X$ onto $\im(\delta^1)$. By \cite{Bo87},
Section II.1, Corollary 1 to Proposition 5, there is a finite subset
$T'\subeq T$ such that
\begin{align*}
\|[c]\|'
\leq& C\sup_{t\in T'}p_t(\widetilde{\delta^1}([c]))
=C\sup_{t\in T'}p_t(\delta^1(c))
=C\sup_{t\in
T'}\|\delta^1_{T'}(c)(t)\|
=C\|\delta^1_{T'}(c)\|_{T'},
\end{align*}
for all $c\in\Map(S,V)$, where $C>0$ is a constant. If
$c\in\ker(\delta^1_{T'})$, the above inequality implies that
$\|[c]\|'=0$. So $c\in\ker(\delta^1)$. Then
$\ker(\delta^1_{T'})\subeq\ker(\delta^1)$. The converse inclusion is
obvious.

Let $\widetilde{\delta^1_{T'}}:X\rightarrow\Map(T',V)$ be the
continuous linear map induced by $\delta^1_{T'}$. The above
inequality says nothing but $\|[c]\|'\leq
C\|\widetilde{\delta^1_{T'}}([c])\|_{T'}$ for $[c]\in X$. So
$\im(\delta^1_{T'})=\im(\widetilde{\delta^1_{T'}})$ is closed in
$\Map(T',V)$. This proves the lemma.
\end{proof}

\textbf{Proof of Theorem~\ref{thm:rigidity}.} Since the Lie algebra
$\L(G)$ of $G$ is a continuous Banach $\Gamma$-module, and we have
assumed that $\delta^1$ has closed image, by Lemma~\ref{lem:finite},
there is a finite subset $T'\subeq T$ such that the map
$\delta^1_{T'}:\Map(S,\L(G))\rightarrow\Map(T',\L(G))$ defined by
$\delta^1_{T'}(c)=\delta^1(c)|_{T'}$ satisfies
$\ker(\delta^1_{T'})=\ker(\delta^1)$, and $\im(\delta^1_{T'})$ is
closed in $(\Map(T',\L(G)),\|\cdot\|_{T'})$.

For a finite set $A$, the set $\Map(A,G)$ of all maps
from $A$ to $G$ is isomorphic to the product $G^A$ of
$|A|$ copies of $G$, hence carries a natural Banach manifold
structure, and the tangent space of $\Map(A,G)$ at the constant map
$e_A:a\mapsto e$ may be identified with $\Map(A,\L(G))$.

Define the maps $\varphi:G\rightarrow\Map(S,G)$ and
$\psi:\Map(S,G)\rightarrow\Map(T',G)$ by
$$\varphi(g)(s)=gr(\pi(s))g^{-1}r(\pi(s))^{-1},$$
$$\psi(\alpha)(t)=(\alpha(s_{t,1})r(\pi(s_{t,1})))^{\epsilon_{t,1}}\cdots
(\alpha(s_{t,m_t})r(\pi(s_{t,m_t})))^{\epsilon_{t,m_t}}.$$ Then it
is easy to see that $\varphi(e)=e_S$, $\psi(e_S)=e_{T'}$, and that
$\psi\circ\varphi\equiv e_{T'}$. Now we compute the differentials of
$\varphi$ and $\psi$ at $e$ and $e_S$, respectively. We have
\begin{align*}
&d\varphi(e)(v)(s)
=\D\varphi(\exp(\lambda v))(s)\\
=&\D\exp(\lambda v)r(\pi(s))\exp(-\lambda v)r(\pi(s))^{-1}
=v-\Ad(r(\pi(s)))v=\delta^0(v)(s),
\end{align*}
where $v\in\L(G)$, $s\in S$. So $$d\varphi(e)=\delta^0.$$ We also
have
\begin{align*}
&d\psi(e_S)(c)(t)
=\D\psi(\exp(\lambda c))(t)\\
=&\D(\exp(\lambda c(s_{t,1}))r(\pi(s_{t,1})))^{\epsilon_{t,1}}\cdots
(\exp(\lambda c(s_{t,m_t}))r(\pi(s_{t,m_t})))^{\epsilon_{t,m_t}}\\
=&\D\prod_{j=1}^{m_t} r(\pi(s_{t,1}^{\epsilon_{t,1}}\cdots
s_{t,j-1}^{\epsilon_{t,j-1}}s_{t,j}^{\epsilon'_{t,j}}))\exp(\epsilon_{t,j}\lambda
c(s_{t,j})) r(\pi(s_{t,1}^{\epsilon_{t,1}}\cdots
s_{t,j-1}^{\epsilon_{t,j-1}}s_{t,j}^{\epsilon'_{t,j}}))^{-1}\\
=&\sum_{j=1}^{m_t}\epsilon_{t,j}\Ad(r(\pi(s_{t,1}^{\epsilon_{t,1}}\cdots
s_{t,j-1}^{\epsilon_{t,j-1}}s_{t,j}^{\epsilon'_{t,j}})))c(s_{t,j})\\
=&\delta^1(c)(t),
\end{align*}
where $c\in\Map(S,\L(G))$, $t\in T'$. So
$d\psi(e_S)=\delta^1_{T'}.$ Since $H^1(\Gamma,\L(G))=\{0\}$, we have
$$\im(d\varphi(e))=\im(\delta^0)=\ker(\delta^1)=\ker(\delta^1_{T'})=\ker(d\psi(e_S)).$$
We also have that $\im(d\psi(e_S))=\im(\delta^1_{T'})$ is closed in
$\Map(T',\L(G))$. This verifies all the conditions of Theorem
\ref{thm:implicit}. So by Theorem~\ref{thm:implicit}, there exists a
neighborhood $W\subeq\Map(S,G)$ of $e_S$ such that
$$\psi^{-1}(e_{T'})\cap W=\varphi(G)\cap W.$$

Let $U$ be the neighborhood of $r$ in $\Hom(\Gamma,G)$ such that
$r'\in U$ if and only if the map $\alpha_{r'}:S\rightarrow G$
defined by $\alpha_{r'}(s)=r'(\pi(s))r(\pi(s))^{-1}$ is in $W$. Then
for any $r'\in U$, $\alpha_{r'}\in\psi^{-1}(e_{T'})\cap W$. By what
we have proved above, $\alpha_{r'}\in\varphi(G)$, which means that
there exists $g\in G$ such that
$$gr(\pi(s))g^{-1}r(\pi(s))^{-1}=r'(\pi(s))r(\pi(s))^{-1}$$ for all
$s\in S$, that is, $r'(\gamma)=gr(\gamma)g^{-1}$ for
$\gamma\in\pi(S)$. But $\pi(S)$ generates $\Gamma$. So
$r'(\gamma)=gr(\gamma)g^{-1}$ for any $\gamma\in\Gamma$. This shows
that $r$ is locally rigid. \qed

\begin{remark}\label{r:hausdorff}
The condition that $\delta^1$ has closed image in Theorem
\ref{thm:rigidity} can be replaced by a more intrinsic condition. We
recall the definition of group cohomology by projective resolutions.
Let $\Gamma$ be a group, let
$$\cdots\rightarrow
P_2\stackrel{\partial_1}\rightarrow
P_1\stackrel{\partial_0}\rightarrow
P_0\stackrel{\varepsilon}\rightarrow \Z\rightarrow0$$ be a
projective $\Gamma$-resolution of the trivial $\Gamma$-module $\Z$.
Then for a $\Gamma$-module $V$, the cohomology $H^*(\Gamma,V)$ of
$\Gamma$ in $V$ is, by definition, the cohomology of the complex
$$0\rightarrow\Hom_\Gamma(P_0,V)\stackrel{\delta^0}\rightarrow\Hom_\Gamma(P_1,V)
\stackrel{\delta^1}\rightarrow\Hom_\Gamma(P_2,V)\rightarrow\cdots,$$
and $H^*(G,V)$ is independent of the particular choice of the
projective resolution. The most usual definition of group cohomology
that appeared in Remark~\ref{r:coincide} corresponds to the standard
resolution. If each $P_j$ in the resolution is countable, and $V$ is
a continuous Banach $\Gamma$-module, then each $\Hom_\Gamma(P_j,V)$,
viewed as a closed subspace of $\Map(P_j,V)$, is a Fr\'{e}chet
space, and the coboundary operators $\delta^j$ are continuous. For
some $n\geq0$, to say that $\delta^n$ has closed image is equivalent
to say that $H^{n+1}(G,V)$ is Hausdorff (see \cite{BW}). Using a
standard argument of homological algebra, it can be shown that
$H^{n+1}(G,V)$ being Hausdorff is independent of the particular
choice of the countable projective resolution. The operator
$\delta^1$ that we used in Theorem~\ref{thm:rigidity} is just the
first coboundary operator defined by the Gruenberg resolution of
$\Gamma$ associated with the presentation of $\Gamma$ (see
\cite{Gr60}), and the modules $P_j$ in the Gruenberg resolution are
countable if the free group $F$ in the presentation of $\Gamma$ is
countable. So the condition that $\delta^1$ has closed image in
Theorem~\ref{thm:rigidity} is equivalent to the condition that
$H^2(\Gamma,\L(G))$ is Hausdorff, which is independent of the
particular countable resolution that we use in the definition of the
group cohomology. By \cite{BW}, Chapter IX, Proposition 3.5, if
$H^2(\Gamma,\L(G))$ is finite-dimensional, then $H^2(\Gamma,\L(G))$
is Hausdorff.
\end{remark}

In view of the above remark, the following corollary of
Theorem~\ref{thm:rigidity} is obvious.

\begin{corollary}
Let $G$ be a Banach--Lie group, $\Gamma$ be a finitely generated
group and $r\in\Hom(\Gamma,G)$. If $H^1(\Gamma,\L(G))=\{0\}$ and
$H^2(\Gamma,\L(G))$ is finite-dimensional, then $r$ is locally
rigid.
\end{corollary}


\begin{thebibliography}{aaaaaa}

\bibitem[AW05]{AW05} An, J. and Wang, Z., {\it Nonabelian cohomology with coefficients in Lie
groups},  Trans. Amer. Math. Soc., to appear.

\bibitem[Be00]{Be00} Benveniste, E. J., {\it Rigidity of isometric lattice actions on
compact Riemannian manifolds}, Geom. Funct. Anal. {\bf 10:3} (2000), 516--542.

\bibitem[Bou87]{Bo87} Bourbaki, N., ``Topological vector spaces, Chapters 1--5," Springer-Verlag, Berlin, 1987.

\bibitem[Bou89]{Bou89} Bourbaki, N., ``Lie Groups and Lie Algebras, Chapter 1--3,''
Springer-Verlag, Berlin, 1989

\bibitem[Bro76]{Bro76} Browder, F. E.,  ``Nonlinear Operators and Nonlinear
Equations of Evolution in Banach Spaces,''  Proc. Symp. Pure Math.
{\bf XVIII:2}, Amer. Math. Soc., 1976.

\bibitem[BW]{BW} Borel, A., Wallach, N., ``Continuous cohomology, discrete subgroups,
and representations of reductive groups", Second edition, American
Mathematical Society, Providence, RI, 2000.


\bibitem[Fi05]{Fi05} Fisher, D., {\it First cohomology and local rigidity of group
actions}, preprint, math.DG/0505520.

\bibitem[GN03]{GN03} Gl\"ockner, H., and K.-H. Neeb, {\it Banach--Lie quotients, enlargibility,
and universal complexifications}, J. Reine Angew. Math. {\bf 560} (2003), 1--28


\bibitem[Gr60]{Gr60} Gruenberg, K. W., {\it Resolutions by relations}, J. London Math. Soc., {\bf
35} (1960), 481--494.

\bibitem[Ha77]{Ha77} Hamilton, R. S., {\it Deformation of complex structures on manifolds
with boundary I: The stable case}, J. Differential Geometry {\bf 12:1} (1977), 1--45.

\bibitem[Ha82]{Ha82} Hamilton, R. S., {\it The inverse function theorem of Nash and
Moser}, Bull. Amer. Math. Soc. (N.S.) {\bf 7:1} (1982), 65--222.

\bibitem[Hof75]{Hof75} Hofmann, K. H.,
{\it Analytic groups without analysis}, Symposia Mathematica {\bf 16}
(Convegno sui Gruppi Topologici e Gruppi di Lie, INDAM, Rome, 1974); 357--374,
Academic Press, London, 1975
%

\bibitem[La99]{La99}
S. Lang, ``Fundamentals of Differential Geometry,''
Graduate Texts in Math. \textbf{191},
Springer-Verlag, Berlin, 1999.

\bibitem[Ne04]{Ne04} Neeb, K.-H., {\it Infinite-dimensional Lie groups and their
representations}, in ``Lie Theory: Lie Algebras and Representations,''
Progress in Math. {\bf 228}, Ed. J.~P.~Anker, B.~\O{}rsted,
Birkh\"auser Verlag, 2004; 213--328

\bibitem[Ne06]{Ne06} K. H. Neeb, {\em Towards a Lie theory of locally convex
groups}, Jap. J. Math. 3rd series {\bf 1:2} (2006), 291--468.

\bibitem[Omo74]{Omo74} H. Omori, ``Infinite-Dimensional Lie
Transformations Groups,''
Lecture Notes in Math. {\bf 427}, Springer-Verlag, 1974.

\bibitem[Ru73]{Ru73} Rudin, W., ``Functional Analysis,'' McGraw Hill, 1973.

\bibitem[We64]{We64} Weil, A., {\it Remarks on the cohomology of groups}, Ann. of Math.
(2), {\bf 80} (1964), 149--157.



\end{thebibliography}
\end{document}